\documentclass[reqno,11pt]{amsart}
\usepackage{dualizing}

\usepackage[matrix,arrow,line,curve,graph,frame]{xy}
\IfFileExists{xypdf.sty}{\usepackage{xypdf}}{}

\newcommand{\CartS}[1][dr]{\ar@{}[#1]|-*+{\square}}

\begin{document}

\title{The $\ell$-adic Dualizing Complex on an Excellent Surface with Rational Singularities}
\author{Ting Li}
\address{Department of Mathematics\\
	Sichuan University\\
	Chengdu 610064\\
	P. R. China}
\email{moduli@live.cn}

\maketitle

\begin{abstract}
In this article, we show that if $X$ is an excellent surface with rational singularities, the constant sheaf $\mathbb{Q}_{\ell}$ is a dualizing complex. In coefficient $\mathbb{Z}_{\ell}$, we also prove that the obstruction for $\mathbb{Z}_{\ell}$ to become a dualizing complex lying on the divisor class groups at singular points. As applications, we study the perverse sheaves and the weights of $\ell$-adic cohomology groups on such surfaces.
\end{abstract}
\plainfoot{Mathematics Subject Classification 2000: 14F20, 14J17}
\plainfoot{This work is supported by NSFC(10626036)}

\section*{Introduction}

In \cite[I]{SGA5}, the dualizing complexes in \'etale cohomology was considered by Grothendieck. O. Gabber made a breakthrough on this fields. He prove that every excellent schemes admits a dualizing complex; and on a regular excellent scheme, the constant sheaf $\mathbb{Z}_{\ell}$ is a dualizing complex. See \cite{JRiou1} for the detail.

In this paper, we study the properties of dualizing complexes on surfaces with rational singularities. In \cite{MR0276239}, J. Lipman prove that a two-dimensional normal local ring $R$ has a rational singularity if and only if $R$ has a finite divisor class group. This makes me think that rational singularities only affect the torsion part of the \'etale cohomology with coefficient $\mathbb{Z}_{\ell}$, but leave the free part invariant. As an evidence, we prove that on an excellent surface $X$ with at most rational singularities, the potential dualizing complex $\mathscr{K}_X$ (determinated by the dimension function $x \mapsto 2 - \dim \mathcal{O}_{X,x}$) concentrates on the entries $-2$ and $-4$; in detail, $\upH{-4}(\mathscr{K}_X) = \mathbb{Z}_{\ell}(2)$ and $\upH{-2}(\mathscr{K}_X)(-1)$ is the $\ell$-torsion part of the divisor class groups at singular points (see \autoref{Main-P0}). In particular, $\mathbb{Q}_{\ell}$ is a dualizing complex of such surface in coefficient $\mathbb{Q}_{\ell}$.

The paper is organized as follows. In \S1, we briefly review the theory of ``potential dualizing complex'' in \cite{JRiou1}. Because the dualizing complexes on a scheme are not unique, they vary by tensor products with invertible complexes. So we prefer the \emph{potential dualizing complex} which is a dualizing complex unique determinated by the dimension function. Basing on this theory, we reiterate the \'etale homology in \cite{MR0444667}. Since lacking of suitable references related to \'etale homology both in adic coefficients and over arbitrary base schemes, it is necessary to spend a little more space on rewriting it. In \S2, we calculate the \'etale homology on curves.

In \S3, we prove the main results. First we calculate the \'etale homology on arbitrary surfaces. Basing on it, we obtain the main results about the dualizing complexes on an excellent surface with rational singularities; and prove the Poincar\'e Duality for such surfaces. In \S4, we give two applications of above theory. First we study perverse sheaves on a surface with rational singularities. Basing on it, we prove that the cohomologies of smooth sheaves of  punctually pure weights on such surfaces, are also punctually pure of weights.

\section{Review of dualizing sheaves and \'etale homology}

Let $\ell$ be a prime number, $\Lambda$ the integral closure of $\mathbb{Z}_{\ell}$ in a finite extension field of $\mathbb{Q}_{\ell}$.

We consider the following conditions related to a scheme $X$:
\begin{itemize}
	\item[(\dag)] $X$ is Noetherian, excellent, of finite Krull dimension, $\ell$ is invertible on $X$ and $\mathrm{cd}_{\ell}(X) < \infty$.
\end{itemize}

From the Gabber's finitenes theorem for \'etale cohomology in \cite{OGab1}, we may construct an adic formalism $\conDb{\etSite{X}, \Lambda}$ together with Grothendieck's six operations on schemes satisfying (\dag). See also \cite{MR1106899} for details. Moreover if $X$ is a scheme satisfying (\dag), then any scheme of finite type over $X$ satisfies (\dag).

Next, we introduces the ``potential dualizing complex'' defined in \cite{JRiou1}.

\begin{Defi}[\cite{JRiou1}, 2.2]
	Let $X$ be a scheme satisfying (\dag), $\delta \colon X \to \mathbb{Z}$ a dimension function. A \emph{potential dualizing complex} (in coefficient $\Lambda$) on $X$ consisting of
	\begin{enumerate}
		\item an object $\mathscr{K}$ in $\conDp{\etSite{X}, \Lambda}$,
		\item an isomorphism $\rdF \varGamma_x(\mathscr{K}) \isoTo \Lambda\bigl(\delta(x)\bigr)\bigl[2\delta(x)\bigr]$ in $\conDp{\etSite{x}, \Lambda}$ for every point $x$ on $X$;
	\end{enumerate}
	and these data satisfy that: for every immediate specialization $\bar{y} \to \bar{x}$ on $X$, the following diagram commutes,
	\[ \xymatrix@C+2em{\rdF \varGamma_{\bar{y}}(\mathscr{K}) \ar[r]^-{\mathrm{sp}^X_{\bar{y} \to \bar{x}}} \ar[d]_-{\cong}
			& \rdF \varGamma_{\bar{x}}(\mathscr{K})(1)[2] \ar[d]^-{\cong} \\
		R\bigl(\delta(y)\bigr)\bigl[2\delta(y)\bigr] \ar@{=}[r] & R\bigl(\delta(x)+1\bigr)\bigl[2(\delta(x)+1)\bigr]} \]
	where $\mathrm{sp}^X_{\bar{y} \to \bar{x}}$ is the transition morphism of codimension $1$ defined in \cite[\S1]{JRiou1}.
\end{Defi}

The potential dualizing complexes have the following properties.

\begin{Prop}[\cite{JRiou1}, 4.1 \& 5.1]
	Every scheme $X$ satisfying (\dag) equipped with a dimension function $\delta$ has a potential dualizing complex, unique up to unique isomorphism, which we denoted by $\mathscr{K}_{X, \delta}$ or simply by $\mathscr{K}_X$. Moreover the potential dualizing complex is a dualizing complex.
\end{Prop}

\begin{Prop}[\cite{JRiou1}, 2.8]\label{236P2}
	Let $X$ be a regular scheme satisfying (\dag). Then $\delta(x) \coloneqq - \dim \mathcal{O}_{X,x}$ is a dimension function on $X$ and we have $\mathscr{K}_{X, \delta} = \Lambda$.
\end{Prop}

\begin{Prop}[\cite{JRiou1}, 4.3]\label{P-78E}
	Let $f \colon X \to Y$ be a compactifiable morphism of schemes satisfying (\dag). Let $\delta_Y$ be a dimension function on $Y$ and equip $X$ with the dimension function
	\[ \delta_X(x) \coloneqq \delta_Y\bigl(f(x)\bigr) + \trD k(x) \big/ k\bigl(f(x)\bigr) \,. \]
	Then we have $\mathscr{K}_{X, \delta_X} = \rdF f^! \mathscr{K}_{Y, \delta_Y}$.
\end{Prop}

\begin{Prop}[\cite{JRiou1}, 4.2]
	Let $f \colon X \to Y$ be a regular morphism of schemes satisfying (\dag). Let $\delta_Y$ be a dimension function on $Y$ and equip $X$ with the dimension function
	\[ \delta_X(x) \coloneqq \delta_Y(y) - \dim \mathcal{O}_{X_y, x} \qquad \bigl(y \coloneqq f(x)\bigr) \]
	Then we have $\mathscr{K}_{X, \delta_X} = f^{\ast} \mathscr{K}_{Y, \delta_Y}$.
\end{Prop}

Basing on above theory, we may generalize the \'etale homology in \cite{MR0444667} to arbitrary schemes satisfying (\dag). Since most proof are almost the same with \cite{MR0444667}, we only give these which meed special care.

\begin{Defi}
	Let $X$ be a scheme satisfying (\dag) equipped with a dimension function $\delta$. For each $n \in \mathbb{Z}$, we define
	\[ \rmH{n}(X, \delta, \Lambda) \coloneqq \upH{-n}(\etSite{X}, \mathscr{K}_{X, \delta}) = \Hom_{\conDb{\etSite{X}, \Lambda}}\bigl(\Lambda_X, \mathscr{K}_{X, \delta}[-n]\bigr) \,. \]

	We also use $\rmH{n}(X, \delta)$ or $\rmH{n}(X, \Lambda)$ or $\rmH{n}(X)$ to denote $\rmH{n}(X, \delta, \Lambda)$, if no confusion arise.
\end{Defi}

Let $f \colon X \to Y$ be a compactifiable morphism of schemes satisfying (\dag), $\delta$ a dimension function on $Y$, and $\delta'$ the dimension function on $X$ induced by $\delta$ as in \autoref{P-78E}. Then we also use $\rmH{n}(X, \delta, \Lambda)$ or $\rmH{n}(X, \delta)$ to denote $\rmH{n}(X, \delta', \Lambda)$.

\medskip

If $f \colon X \to Y$ is a proper morphism of schemes satisfying (\dag), $\delta$ a dimension function on $Y$. For each $n \in \mathbb{Z}$, we define a homomorphism of $\Lambda$-modules
\[ f_{\ast} \colon \rmH{n}(X, \delta) \to \rmH{n}(Y, \delta) \]
as follows. For each $\alpha \in \rmH{n}(X, \delta)$, regarding $\alpha \colon \Lambda_X \to \mathscr{K}_{X, \delta}[-n]$ as a morphism in $\conDb{\etSite{X}, \Lambda}$, then $f_{\ast}(\alpha)$ is defined to be the composition
\[ \Lambda_Y \to \rdF f_{\ast} \Lambda_X \xrightarrow{\rdF f_{\ast}(\alpha)} \rdF f_{\ast} \mathscr{K}_{X, \delta}[-n] \isoTo \rdF f_{\ast} \circ f^{\ast} \mathscr{K}_{Y, \delta}[-n] \to \mathscr{K}_{Y, \delta}[-n] \,, \]
where the first and the last morphisms are induced by the adjunctions.

It is easy to verify that if $f \colon X \to Y$ and $g \colon Y \to Z$ are two proper morphisms, then $(g \circ f)_{\ast} = g_{\ast} \circ f_{\ast}$.

If $Y$ is the spectrum of a separably closed field, as $\rmH{0}(Y) = \Lambda$, we may write
\[ \deg \coloneqq f_{\ast} \colon \rmH{0}(X) \to \Lambda \,. \]

\begin{Prop}\label{T-P0}
	Let $X$ be a scheme satisfying (\dag) equipped with a dimension function $\delta$, $Y$ a closed subscheme of $X$ and $U \coloneqq X \setminus Y$. Then there is a long exact sequence of $\Lambda$-modules:
	\[ \cdots \rmH{n+1}(U, \delta) \to \rmH{n}(Y, \delta) \to \rmH{n}(X, \delta) \to \rmH{n}(U, \delta) \to \rmH{n-1}(Y, \delta) \to \cdots \,. \]
\end{Prop}

\begin{Prop}[Mayer-Vietoris Sequence]\label{T-P1}
	Let $X$ be a scheme satisfying (\dag) equipped with a dimension function $\delta$, $X_1$ and $X_2$ two closed subschemes of $X$ such that $X = X_1 \cup X_2$ (as sets). Then we have a long exact sequence
	\[ \cdots \to \rmH{n}(X_1 \cap X_2, \, \delta) \to \rmH{n}(X_1, \delta) \oplus \rmH{n}(X_2, \delta) \to \rmH{n}(X, \delta) \to \rmH{n-1}(X_1 \cap X_2, \, \delta) \to \cdots \,. \]
\end{Prop}

\begin{Nota}
	Let $X$ be a Noetherian scheme, $\delta$ a dimension function on $X$. Then we define
	\[ \dim_{\delta}(X) \coloneqq \sup_{x \in X}\delta(x) \,. \]
	If $\xi_1, \xi_2, \ldots, \xi_r$ are all generic points of irreducible components of $X$, then
	\[ \dim_{\delta}(X) = \max_{1 \leqslant i \leqslant r} \delta(\xi_i) < +\infty \,. \]
\end{Nota}

\begin{Prop}[Vanishing]
	Let $X$ be a scheme satisfying (\dag) equipped with a dimension function $\delta$. Then $\rmH{n}(X, \delta) = 0$ for all $n > 2 \dim_{\delta}(X)$.
\end{Prop}

\begin{Cor}\label{53P2}
	Let $X$ be a scheme satisfying (\dag) equipped with a dimension function $\delta$, $Y$ a closed subscheme of $X$, $X' \coloneqq X \setminus Y$. Then for each integer $n > 2 \dim_{\delta}(Y) + 1$, there is a canonical isomorphism $\rmH{n}(X, \delta) \isoTo \rmH{n}(X', \delta)$ of $\Lambda$-modules.
\end{Cor}

\begin{proof}
	We have only to apply \autoref{T-P0}.
\end{proof}

We define the pull-back maps as follows. Let $Y$ be a scheme satisfying (\dag) equipped with a dimension function $\delta$, $f \colon X \to Y$ is a compactifiable flat morphism of relative dimension $d$. For each $n \in \mathbb{Z}$, we define a homomorphism of $\Lambda$-modules
\[ f^{\ast} \colon \rmH{n}(Y, \delta) \to \rmH{n+2d}(X, \delta)(-d) \]
as follows. For each $\beta \in \rmH{n}(Y, \delta)$, $f^{\ast}(\beta)$ is defined to be the composition
\[ \Lambda_X \xrightarrow{\mathrm{t}_f} \rdF f^! \Lambda_X[-2d](-d) \xrightarrow{\rdF f^!(\beta)} \rdF f^! \mathscr{K}_Y[-n-2d](-d) \isoTo \mathscr{K}_X[-n-2d](-d) \,, \]
where $\mathrm{t}_f$ is the morphism dual to the trace morphism $\mathrm{Tr}_f$ (see \cite[XVIII (3.2.1.2)]{SGA4}).

It is easy to verify that if $f \colon X \to Y$ and $g \colon Y \to Z$ be two compactifiable flat morphisms of relative dimension $d$ and $e$ respectively, then
\[ (g \circ f)^{\ast} = f^{\ast} \circ g^{\ast} \colon \rmH{n}(Z) \to \rmH{n+2d+2e}(X)(-d-e) \,. \]

Moreover the maps $f_{\ast}$ and $f^{\ast}$ commute in Cartesian squares.

Now we define the cycle map for \'etale homology.

\begin{Nota}
	Let $X$ be a scheme satisfying (\dag) equipped with a dimension function $\delta$, $d \coloneqq \dim_{\delta}(X)$.
	\begin{enumerate}
		\item First we assume that $X$ is regular and integral. Then we use $\mathrm{t}_X \colon \Lambda \isoTo \mathscr{K}_{X, \delta}(-d)[-2d]$ to denote the canonical morphism defined by \autoref{236P2}, and use $\clM(X)$ to denote the element in $\rmH{2d}(X, \delta)(-d)$ corresponding to $\mathrm{t}_X$.
		\item Second we only assume that $X$ is integral. As $X$ is excellent, $U \coloneqq X_{\mathrm{reg}}$ is an open dense subset of $X$. Obviously $\dim_{\delta}(X \setminus U) < d$. By \autoref{53P2}, there is a canonical isomorphism
			\[ \rmH{2d}(X, \delta)(-d) \isoTo \rmH{2d}(U, \delta)(-d) \,. \]
			Let $\clM(X)$ be the inverse image of $\clM(U)$ under above isomorphism.
		\item In general case, we let $X_1, X_2, \ldots, X_r$ be all irreducible components of $X$ with $\dim_{\delta}(X_i) = d$. For each $i$, regard $X_i$ as a reduced subscheme of $X$, let $\iota_i \colon X_i \hookrightarrow X$ be the inclusion and let $\xi_i$ be the generic point of $X_i$. Then we define
			\[ \clM(X) \coloneqq \sum^r_{i=1} \mathrm{length}(\mathcal{O}_{X, \xi_i}) \cdot \iota_{i, \ast} \bigl(\clM(X_i)\bigr) \in \rmH{2d}(X, \delta)(-d) \,. \]
			We also use $\mathrm{t}_X \colon \Lambda \isoTo \mathscr{K}_{X, \delta}(-d)[-2d]$ to denote the morphism in $\conDb{X, \Lambda}$ corresponding to $\clM(X)$.
	\end{enumerate}
\end{Nota}

\begin{Nota}
	Let $X$ be a scheme satisfying (\dag) equipped with a dimension function $\delta$, $Y$ a closed subscheme of $X$, $d \coloneqq \dim_{\delta}(Y)$. Then we define
	\[ \clM_X(Y) \coloneqq \iota_{\ast}\bigl(\clM(Y)\bigr) \in \rmH{2d}(X, \delta)(-d) \,, \]
	where $\iota \colon Y \hookrightarrow X$ is the inclusion.
\end{Nota}

Next, we study the \'etale homology under birational morphisms. First we need two lemmas.

\begin{Lem}\label{194L0}
	Let $X$ be a scheme satisfying (\dag), $Y$ a closed subscheme of $X$, $U \coloneqq X \setminus Y$, $i \colon Y \hookrightarrow X$ and $j \colon U \hookrightarrow X$ the inclusions. Then for each object $\mathscr{F}$ in $\conDb{\etSite{X}, \Lambda}$, there are two distinguished triangles in $\conDb{\etSite{X}, \Lambda}$:
	\begin{enumerate}
		\item $j_!j^{\ast}\mathscr{F} \to \mathscr{F} \to i_{\ast}i^{\ast}\mathscr{F} \to j_!j^{\ast}\mathscr{F}[1]$\,,
		\item $i_{\ast}\rdF i^!\mathscr{F} \to \mathscr{F} \to \rdF j_{\ast} j^{\ast}\mathscr{F} \to i_{\ast}\rdF i^!\mathscr{F}[1]$\,.
	\end{enumerate}
\end{Lem}

\begin{Lem}
	Let $X$ be a scheme satisfying (\dag), $Y$ a closed subscheme of $X$, $U \coloneqq X \setminus Y$, $i \colon Y \hookrightarrow X$ and $j \colon U \hookrightarrow X$ the inclusions, $\mathscr{F}' \xrightarrow{\varphi} \mathscr{F} \xrightarrow{\psi} \mathscr{F}''$ a sequence of objects in $\conDb{\etSite{X}, \Lambda}$. Assume that one of the following two conditions holds:
	\begin{enumerate}
		\item $j^{\ast}(\varphi) \colon j^{\ast}\mathscr{F}' \isoTo j^{\ast}\mathscr{F}$ is an isomorphism, $j^{\ast}\mathscr{F}'' = 0$, and for all $n \in \mathbb{Z}$,
			\[ 0 \to \upH{n}(i^{\ast}\mathscr{F}') \to \upH{n}(i^{\ast}\mathscr{F}) \to  \upH{n}(i^{\ast}\mathscr{F}'') \to 0 \]
			is a short exact sequence;
		\item $j^{\ast}(\psi) \colon j^{\ast}\mathscr{F} \isoTo j^{\ast}\mathscr{F}''$ is an isomorphism, $j^{\ast}\mathscr{F}' = 0$, and for all $n \in \mathbb{Z}$,
			\[ 0 \to \upH{n}(\rdF i^!\mathscr{F}') \to \upH{n}(\rdF i^!\mathscr{F}) \to  \upH{n}(\rdF i^!\mathscr{F}'') \to 0 \]
			is a short exact sequence.
	\end{enumerate}
	Then there exists a morphism $\delta \colon \mathscr{F}'' \to \mathscr{F}'[1]$ which makes
	\[ \mathscr{F}' \xrightarrow{\varphi} \mathscr{F} \xrightarrow{\psi} \mathscr{F}'' \xrightarrow{\delta} \mathscr{F}'[1] \]
	a distinguished triangle in $\conDb{\etSite{X}, \Lambda}$.
\end{Lem}

\begin{proof}
	The morphism $\varphi$ extends to a distinguished triangle
	\[ \mathscr{F}' \xrightarrow{\varphi} \mathscr{F} \xrightarrow{\phi} \mathscr{G} \xrightarrow{\rho} \mathscr{F}'[1] \]
	in $\conDb{\etSite{X}, \Lambda}$. Note that $0 \to \mathscr{F}'' \xrightarrow{\iDe} \mathscr{F}'' \to 0[1]$ is also a distinguished triangle. Hence there exists a morphism $\alpha \colon \mathscr{G} \to \mathscr{F}''$ such that the triple $(0, \psi, \alpha)$ is a morphism of distinguished triangles, i.e., the following diagram commutes.
	\[ \xymatrix{\mathscr{F}' \ar[r]^-{\varphi} \ar[d]^{0} & \mathscr{F} \ar[r]^-{\phi} \ar[d]^{\psi} & \mathscr{G} \ar[r]^-{\rho} \ar[d]^{\alpha} & \mathscr{F}'[1] \ar[d]^{0} \\
		0 \ar[r] & \mathscr{F}'' \ar[r]^-{\iDe} & \mathscr{F}'' \ar[r] & 0[1]} \]

	Now if (1) holds, then both $j^{\ast}(\alpha)$ and $i^{\ast}(\alpha)$ are isomorphic; and if (2) holds, then both $j^{\ast}(\alpha)$ and $\rdF i^!(\alpha)$ are isomorphic. In either case, by \autoref{194L0}, the morphism $\alpha \colon \mathscr{G} \isoTo \mathscr{F}''$ is an isomorphism. Therefore we may let $\delta \coloneqq \rho \circ \alpha^{-1}$.
\end{proof}

Basing on above lemma, we have the following two propositions.

\begin{Prop}\label{194P32}
	Let
	\[ \xymatrix{Y' \ar@{^(->}[r]^-{i'} \ar[d]_{q} \CartS & X' \ar[d]^{p} \\ Y \ar@{^(->}[r]^-{i} & X} \]
	be a Cartesian square of schemes satisfying (\dag), $r \coloneqq p \circ i' = i \circ q$. Assume that $p$ is proper, $i$ is a closed immersion; and there exists an open subset $U$ of $X$ such that $i(Y) = X \setminus U$ and $p$ induces an isomorphism $p^{-1}(U) \isoTo U$. Then for each object $\mathscr{F}$ in $\conDb{\etSite{X}, \Lambda}$, there are two distinguished triangles:
	\begin{gather*}
		\mathscr{F} \xrightarrow{\delta_i \oplus \delta_p} i_{\ast} \circ i^{\ast} \mathscr{F} \oplus \rdF p_{\ast} \circ p^{\ast} \mathscr{F} \xrightarrow{\delta_q \circ \prM{1} - \delta_{i'} \circ \prM{2}} \rdF r_{\ast} \circ r^{\ast} \mathscr{F} \to \mathscr{F}[1]\,, \\
		\rdF r_{\ast} \circ \rdF r^! \mathscr{F} \xrightarrow{\vartheta_q \oplus \vartheta_{i'}} i_{\ast} \circ \rdF i^! \mathscr{F} \oplus \rdF p_{\ast} \circ \rdF p^! \mathscr{F} \xrightarrow{-\vartheta_i \circ \prM{1} + \vartheta_p \circ \prM{2}} \mathscr{F} \to (\rdF r_{\ast} \circ \rdF r^! \mathscr{F})[1] \,,
	\end{gather*}
	in $\conDp{\etSite{X}, \Lambda}$, where $\delta$ and $\vartheta$ are induced by the adjunctions.
\end{Prop}

Apply the second triangle in above proposition to $\mathscr{K}_X$, we obtain the following proposition.

\begin{Prop}\label{53P5}
	Let the assumptions and the notations be as in \autoref{194P32}; and let $\delta$ be a dimension function on $X$. Then there is a long exact sequence of $\Lambda$-modules:
	\[ \cdots \to \rmH{n+1}(X, \delta) \to \rmH{n}(Y', \delta) \xrightarrow{q_{\ast} \oplus i'_{\ast}} \rmH{n}(Y, \delta) \oplus \rmH{n}(X', \delta) \xrightarrow{-i_{\ast} \circ \prM{1} + p_{\ast} \circ \prM{2}} \rmH{n}(X, \delta) \to \cdots \,. \]
\end{Prop}

\section{Calculation of \'etale homology on the exceptional divisor}

In this section, we calculate the \'etale homology of exceptional curves.

\begin{Prop}\label{T-P2}
	Let $C$ be a $1$-equidimensional proper algebraic scheme over a separably closed field $k$. Let $C^{(1)}, C^{(2)}, \ldots, C^{(r)}$ be all connected components of $C$; and let $C_1, C_2, \ldots, C_n$ be all irreducible components of $C$. Assume that $\upH{1}(C, \mathcal{O}_C) = 0$. Then we have
	\begin{enumerate}
		\item As to homology, we have
			\begin{enumerate}
				\item There is a canonical isomorphism:
					\[ \rmH{0}(C) = \bigoplus^r_{i=1} \rmH{0}(C^{(i)}) \isoTo[\oplus \deg] \Lambda^{\oplus r} \,. \]
				\item $\rmH{1}(C) = 0$.
				\item $\rmH{2}(C)(-1)$ is a free $\Lambda$-module with basis $\clM(C_1),\ldots,\clM(C_n)$.
			\end{enumerate}
		\item As to cohomology, we have
			\begin{enumerate}
				\item There is a canonical isomorphism:
					\[ \upH{0}(C, \Lambda) = \bigoplus^r_{i=1} \upH{0}(C^{(i)}) \isoTo \Lambda^{\oplus r} \,. \]
				\item $\upH{1}(C, \Lambda) = 0$.
				\item There is a canonical isomorphism $\upH{2}\bigl(C, \Lambda(1)\bigr) \isoTo \Lambda^{\oplus n}$ which sends $c_1(\mathscr{L})$ to $\bigl(\deg(\mathscr{L}|_{C_1}), \ldots, \deg(\mathscr{L}|_{C_r})\bigr)$ for each invertible $\mathcal{O}_C$-module $\mathscr{L}$.
			\end{enumerate}
	\end{enumerate}
\end{Prop}

\begin{proof}
	(1). Let $\bar{k}$ be the algebraic closure of $k$. After replacing $k$ with $\bar{k}$ and $X$ with $(X \otimes_k \bar{k})_{\mathrm{red}}$, we may assume that $k$ is algebraically closed and $X$ is reduced. Now we use induction on the number $n$ of irreducible components of $C$.

	First consider the case $n = 1$. In this case, $C$ is integral. Let $\widetilde{C}$ be the normalization of $C$. Then there is a long exact sequence
	\[ 0 \to \varGamma(C, \mathcal{O}_C) \to \varGamma(C, \mathcal{O}_{\widetilde{C}}) \to \varGamma(C, \mathcal{O}_{\widetilde{C}}/\mathcal{O}_C) \to \upH{1}(C, \mathcal{O}_C) \,. \]
	Note that $\upH{1}(C, \mathcal{O}_C) = 0$, $\varGamma(C, \mathcal{O}_C) = k$ and $\varGamma(C, \mathcal{O}_{\widetilde{C}}) = \varGamma(\widetilde{C}, \mathcal{O}_{\widetilde{C}}) = k$. Hence $\varGamma(C, \mathcal{O}_{\widetilde{C}}/\mathcal{O}_C) = 0$. Since $\mathcal{O}_{\widetilde{C}}/\mathcal{O}_C$ is a coherent $\mathcal{O}_C$-module whose support is of dimension $\leqslant 0$, we have $\mathcal{O}_{\widetilde{C}}/\mathcal{O}_C = 0$, i.e., $C = \widetilde{C}$ is normal. So we have $C = \mathrm{P}_k^1$.

	Assume that (1) is valid for all integers $< n$. Let $C'$ be an irreducible component $C$ and let $C''$ be the union of all other irreducible components. Regard $C'$ and $C''$ as reduced closed subschemes of $C$. Then there is a short exact sequence:
	\[ 0 \to \mathcal{O}_C \xrightarrow{c \mapsto (\bar{c}, \bar{c})} \mathcal{O}_{C'} \oplus \mathcal{O}_{C''} \xrightarrow{(a, b) \mapsto \bar{a}-\bar{b}} \mathcal{O}_{C' \cap C''} \to 0 \,, \]
	which induces a long exact sequence:
	\begin{gather*}
		0 \to \varGamma(C, \mathcal{O}_C) \to \varGamma(C, \mathcal{O}_{C'}) \oplus \varGamma(C, \mathcal{O}_{C''}) \to \varGamma(C, \mathcal{O}_{C' \cap C'}) \\
		\to \upH{1}(C, \mathcal{O}_C) \to \upH{1}(C, \mathcal{O}_{C'}) \oplus \upH{1}(C, \mathcal{O}_{C''}) \to 0 \,.
	\end{gather*}
	Let $r$ and $s$ be the numbers of connected components of $C$ and $C''$ respectively; and let $Z_1, Z_2, \ldots, Z_m$ be all connected components of $C''$ which intersect with $C'$. Then $m =s - r+ 1$. On the other hand, we have
	\[ \dim_k \varGamma(C, \mathcal{O}_C) = r, \quad \dim_k \varGamma(C, \mathcal{O}_{C'}) = 1, \quad \dim_k \varGamma(C, \mathcal{O}_{C''}) = s. \]
	As $\upH{1}(C, \mathcal{O}_C) = 0$, we have $\upH{1}(C, \mathcal{O}_{C'}) = \upH{1}(C, \mathcal{O}_{C''}) = 0$, and
	\[ \dim_k \varGamma(C, \mathcal{O}_{C' \cap C'}) = s-r+1 = m \,. \]
	This shows that for each $i$, $Z_i \cap C'$ contains only one point $P_i$, and $C' \cap C'' = \{P_1, P_2, \ldots, P_m\}$.

	Using induction on $Z_i$, $C'$ and $C''$, we obtain that (1) is valid on all $Z_i$, $C'$ and $C''$. Applying \autoref{T-P1}, we obtain an isomorphism $\rmH{2}(C') \oplus \rmH{2}(C'') \isoTo \rmH{2}(C)$ and an exact sequence
	\[ 0 \to \rmH{1}(C) \to \rmH{0}(C' \cap C'') \xrightarrow{l} \rmH{0}(C') \oplus \rmH{0}(C'') \to \rmH{0}(C) \to 0 \,. \]
	Note that the map $l$ has an inverse which is defined by the composition
	\[ \rmH{0}(C') \oplus \rmH{0}(C'') \xrightarrow{\text{pr}} \rmH{0}(C'') \xrightarrow{\text{pr}} \bigoplus^m_{i=1} \rmH{0}(Z_i) \isoTo[\oplus \varphi_i] \bigoplus^m_{i=1} \rmH{0}(P_i) = \rmH{0}(C' \cap C'') \,, \]
	where the isomorphism $\varphi_i$, which maps $\clM_{Z_i}(P_i)$ to $\clM(P_i)$, is obtained by applying (a) on $Z_i$. Therefore (1) are valid for $C$.

	\medskip

	(2) Let $\pi \colon C \to \Spec k$ be the structural morphism. Note that
	\[ \upH{q}(C, \Lambda) = \upH{q}(\rdF \pi_{\ast} \Lambda_C) = \upH{q}\bigl(\mathscr{D}(\rdF \pi_{\ast} \mathscr{D}(\Lambda_C))\bigr) = \upH{q}\bigl(\mathscr{D}(\rdF \pi_{\ast} \mathscr{K}_C)\bigr) \,, \]
	where $\mathscr{D} \coloneqq \mathscr{H}om(\smdot[.6em], \mathscr{K})$ is the dualizing functor. By (1), we have
	\[ \upH{q}(\rdF \pi_{\ast} \mathscr{K}_C) = \rmH{-q}(C) = \begin{cases} \Lambda^{\oplus r} & q = 0, \\ \Lambda(1)^{\oplus n} & q = -2, \\ 0 & \text{otherwise.} \end{cases} \]
	In other words, all entries of $\rdF \pi_{\ast} \mathscr{K}_C$ are free $\Lambda$-modules. Therefore we have
	\[ \upH{q}(C, \Lambda) = \begin{cases} \Lambda^{\oplus r} & q = 0, \\ \Lambda(-1)^{\oplus n} & q = 2, \\ 0 & \text{otherwise.} \end{cases} \]
	So (2) is verifies. (c) is also by \cite[IX, 4.7]{SGA4}.
\end{proof}

\begin{Rem}
	In fact, every connected component of $(C \otimes_k \bar{k})_{\mathrm{red}}$ is a tree of nonsingular rational curves with normal crossings.
\end{Rem}

\section{Calculation of \'etale homology on surfaces with rational singularities}

In this paper, a surface means a two-dimensional integral normal scheme which satisfies (\dag). Note that for each surface $X$, $\delta(x) \coloneqq 2 - \dim \mathcal{O}_{X, x}$ is a dimension function on $X$; and for all points $x \in X$, we have $\delta(x) = \dim \overline{\{x\}}$. From now on, all surfaces are tacitly equipped with this type of dimension function.

We recall the definition of \emph{rational singularity} (cf. \cite{MR0199191} \& \cite{MR0276239}). A normal local ring $R$ of dimension $2$ is said to have a \emph{rational singularity} if there exists a desingularization $f \colon X \to \Spec(R)$ such that $\upH{1}(X, \mathcal{O}_X) = 0$. By \cite[Proposition (16.5)]{MR0276239}, for a normal local ring $R$ of dimension $2$, if one of the four rings $R$, $\widehat{R}$, $R^{\mathrm{sh}}$ and $\widehat{R^{\mathrm{sh}}}$ has a rational singularity, so have the other three rings. A surface $X$ is said to have a rational singularity at a closed point $x$, if the local ring $\mathcal{O}_{X, x}$ has a rational singularity.

\medskip

Now we calculate the \'etale homology on arbitrary surfaces.

\begin{Prop}\label{T-3P8}
	Let $X$ be a surface which is the spectrum of a strictly Henselian local ring, $P$ the closed point of $X$, $U \coloneqq X \setminus P$. Assume that $U$ is regular, and $X$ admits a desingularization $\pi \colon \widetilde{X} \to X$ which induces an isomorphism $\pi^{-1}(U) \isoTo U$. Let $E \coloneqq \pi^{-1}(P)$ be the exceptional divisor. Then we have
	\begin{enumerate}
		\item The map $\upH{0}(\mathrm{t}_X) \colon \Lambda \isoTo \rmH{4}(X)(-2)$ is an isomorphism.
		\item $\rmH{3}(X) \isoTo \upH{1}\bigl(E, \Lambda(2)\bigr)$.
		\item There is an exact sequence of $\Lambda$-modules:
			\[ 0 \to \rmH{2}(E) \to \upH{2}\bigl(E, \Lambda(2)\bigr) \to \rmH{2}(X) \to \rmH{1}(E) \to 0 \,, \]
			Moreover the kernel of the map $\rmH{2}(X) \to \rmH{1}(E)$ is a finite group.
		\item There is an exact sequence of $\Lambda$-modules:
			\[ 0 \to \rmH{1}(X) \to \rmH{0}(E) \xrightarrow{\deg} \Lambda \to 0 \,. \]
		\item $\rmH{q}(X) = 0$ for $q > 4$ or $q < 1$.
	\end{enumerate}
\end{Prop}

\begin{proof}
	After applying \autoref{53P5} to the following Cartesian square,
	\begin{equation}\label{E-D2}
		\vcenter{\xymatrix{E \ar@{^(->}[r]^-{\tilde{i}} \ar[d]_{\tau} \CartS & \widetilde{X} \ar[d]^{\pi} \\ P \ar@{^(->}[r]^-{i} & X}}
	\end{equation}
	we obtain a long exact sequence
	\begin{equation}\label{E-D3}
		\cdots \to \rmH{n+1}(X) \to \rmH{n}(E) \to \rmH{n}(P) \oplus \rmH{n}(\widetilde{X}) \to \rmH{n}(X) \to \cdots \,.
	\end{equation}
	For each $q \in \mathbb{Z}$, there is a canonical isomorphism of $\Lambda$-modules:
	\[ \rmH{q}(\widetilde{X}) \isoTo \upH{4-q}\bigl(\widetilde{X}, \Lambda(2)\bigr) \isoTo \upH{4-q}\bigl(E, \Lambda(2)\bigr) \,, \]
	where the first isomorphism is induced by \autoref{236P2} (since $\widetilde{X}$ is regular), and the last isomorphism is from the base change theorem via the square \eqref{E-D2}. Thus $\rmH{q}(\widetilde{X}) = 0$ for $q < 2$ and $q > 4$. Also note that $\rmH{q}(E) = 0$ for $q < 0$ and $q > 2$; and $\rmH{q}(P) = \begin{cases} 0 & q \neq 0, \\ \Lambda & q= 0. \end{cases}$ Hence we may divide \eqref{E-D3} to three exact sequences
	\begin{gather} \label{G-E1}
		0 \to \upH{0}\bigl(E, \Lambda(2)\bigr) \to \rmH{4}(X) \to 0, \\ \label{G-E2}
		0 \to \upH{1}\bigl(E, \Lambda(2)\bigr) \to \rmH{3}(X) \to \rmH{2}(E) \to \upH{2}\bigl(E, \Lambda(2)\bigr) \to \rmH{2}(X) \to \rmH{1}(E) \to 0, \\ \label{G-E3}
		0 \to \rmH{1}(X) \to \rmH{0}(E) \xrightarrow{\deg} \Lambda \to \rmH{0}(X) \to 0.
	\end{gather}
	Now (1) is by the exact sequence \eqref{G-E1}.

	From the sequence \eqref{G-E3}, to prove (4) \& (5) we have only to prove the map $\rmH{0}(E) \xrightarrow{\deg} \Lambda$ is epimorphic. Since $X$ is strictly local, the field $k \coloneqq k(P)$ is separably closed. Let $Q$ be a closed point on $E$. Then $\deg \clM_E(Q) = \bigl[k(Q), k\bigr]$. Put $p \coloneqq \Char k$. Since $k$ is separably closed, the number $\bigl[k(Q), k\bigr]$ is either $1$ or a power of $p$; in both cases, it is invertible on $\Lambda$ (as $\ell \neq p$). Hence $\deg$ is surjective.

	From the sequence \eqref{G-E2}, to prove (2) \& (3) we have only to prove that the homomorphism $\rmH{2}(E)(-1) \to \upH{2}\bigl(E, \Lambda(1)\bigr)$ is injective and its cokernel is a torsion group. Let $E_1, E_2, \ldots, E_n$ be all irreducible components of $E$. Then $\rmH{2}(E)(-1)$ is a free $\Lambda$-Module with a basis $\bigl\{\clM_E(E_i)\bigr\}^n_{i=1}$. Note that the composite map
	\begin{equation}\label{E-D1}
		\rmH{2}(E)(-1) \to \rmH{2}(\widetilde{X})(-1) \isoTo \upH{2}\bigl(\widetilde{X}, \Lambda(1)\bigr) \isoTo \upH{2}\bigl(E, \Lambda(1)\bigr) \isoTo[\rho] \Lambda^{\oplus n}
	\end{equation}
	sends each $\clM_E(E_i)$ to $\sum\limits^n_{j=1}a_{ij}e_j$, where the isomorphism $\rho$ is defined in \cite[IX, 4.7]{SGA4}, $e_1, e_2, \ldots, e_n$ is the standard basis for $\Lambda^{\oplus n}$, and
	\[ a_{ij} = \deg \bigl(\mathscr{L}_{\widetilde{X}}(E_i)|_{E_j}\bigr) = (E_i, E_j) \,, \]
	where $\mathscr{L}_{\widetilde{X}}(E_i)$ means the invertible sheaf on $\widetilde{X}$ associated to the divisor $E_i$. In other words, the map \eqref{E-D1} is given by the intersection matrix $\bigl((E_i, E_j)\bigr)$, which is negative-definite by \cite[(14.1)]{MR0276239}. Therefore the map \eqref{E-D1} is injective with torsion cokernel.
\end{proof}

\begin{Thm}\label{Main-P0}
	Let the assumptions and the notations be as in \autoref{T-3P8}. We further assume that $X$ has a rational singularity at $P$. Then we have
	\begin{enumerate}
		\item Let $\mathrm{Cl}(X)$ denote the Weil divisor class group of $X$. Then there is a canonical isomorphism
			\[ \clM \colon \mathrm{Cl}(X) \otimes_{\mathbb{Z}} \Lambda \isoTo \rmH{2}(X)(-1) \]
			which sends each $[Y] \otimes 1$ to $\clM_X(Y)$.
		\item The map $\upH{0}(\mathrm{t}_X) \colon \Lambda \isoTo \rmH{4}(X)(-2)$ is an isomorphism.
		\item $\rmH{q}(X) = 0$ for $q \neq 2, 4$.
	\end{enumerate}
\end{Thm}

\begin{proof}
	By \cite[(4.1)]{MR0276239}, any desingularization of a surface with rational singularities is a product of quadratic transformations. Since $X$ is the spectrum of a local ring, the exceptional divisor $E$ is connected. Thus the map $\deg \colon \rmH{0}(E) \to \Lambda$ is an isomorphism by \autoref{T-P2} (1a). Then (2) \& (3) are by \autoref{T-3P8} and \autoref{T-P2}. So we have only to prove (1).

	We adopt the assumptions and notation in \cite[\S14]{MR0276239}. Let $E_1, E_2, \ldots, E_n$ be all irreducible components of $E$. For each $i$, let $d_i > 0$ be the greatest common divisor of all the degrees of invertible sheaves on $E_i$. Put $k \coloneqq k(P)$ and $p \coloneqq \Char k$. For each $i$, let $Q_i$ be a closed point on $\mathrm{Reg}(E_i)$. Then $Q_i$ defines an invertible sheaf of degree $\bigl[k(Q_i), k\bigr]$ on $E_i$. Since the number $\bigl[k(Q_i), k\bigr]$ is invertible on $\Lambda$, so is $d_i$.

	Let $\mathbf{E}$ be the additive group of divisors on $\widetilde{X}$ generated by $\{E_i\}^n_{i=1}$. Then $\mathbf{E}$ is a free abelian group with basis $E_1, E_2, \ldots, E_n$. Put $\mathbf{E}^{\ast} \coloneqq \Hom_{\mathbb{Z}}(\mathbf{E}, \mathbb{Z})$ and let $\varepsilon_1, \varepsilon_2, \ldots, \varepsilon_n$ be the dual basis. We define a homomorphism of groups:
	\[ \vartheta \colon \mathbf{E} \to \mathbf{E}^{\ast} \,, \qquad E_i \mapsto \sum^n_{j=1} \tfrac{1}{d_j}(E_i, E_j) \cdot \varepsilon_j \,. \]
	By \cite[(14.3), (14.4) \& (17.1)]{MR0276239}, there is an exact sequence
	\[ 0 \to \mathbf{E} \xrightarrow{\vartheta} \mathbf{E}^{\ast} \to \mathrm{Cl}(X) \to 0 \,. \]
	Let $e_1,e_2,\ldots,e_n$ be the basis for the free $\Lambda$-module $\upH{2}\bigl(E, \Lambda(1)\bigr)$ defined by the isomorphism in \autoref{T-P2} (2c). Since each $d_j$ is invertible on $\Lambda$, the homomorphism
	\[ \delta \colon \mathbf{E}^{\ast} \otimes_{\mathbb{Z}} \Lambda \to \upH{2}\bigl(E, \Lambda(1)\bigr) \,, \qquad \varepsilon_j \otimes 1 \mapsto d_j e_j \]
	is an isomorphism of $\Lambda$-modules. It is easy to check that the following diagram
	\[ \xymatrix{0 \ar[r] & \mathbf{E} \otimes_{\mathbb{Z}} \Lambda \ar[r]^{\vartheta \otimes \iDe} \ar[d]^{\text{(i)}}_{\cong}
			& \mathbf{E}^{\ast} \otimes_{\mathbb{Z}} \Lambda \ar[r] \ar[d]^{\delta}_{\cong}
			& \mathrm{Cl}(X) \otimes_{\mathbb{Z}} \Lambda \ar[r] \ar[d]^{\clM} & 0 \\
		0 \ar[r] & \rmH{2}(E)(-1) \ar[r] & \upH{2}\bigl(E, \Lambda(1)\bigr) \ar[r] & \rmH{2}(X)(-1) \ar[r] & 0} \]
	is commutative with both rows exact, where the isomorphism (i) sends each $E_i \otimes 1$ to $\clM_E(E_i)$. Therefore $\clM \colon \mathrm{Cl}(X) \otimes_{\mathbb{Z}} \Lambda \isoTo \rmH{2}(X)(-1)$ is an isomorphism.
\end{proof}

By \cite[(17.1)]{MR0276239}, the divisor class group $\mathrm{Cl}(X)$ is a finite group. So we have:

\begin{Cor}\label{M-A21}
	Let $X$ be a surface having only finitely many singular points, all of which are rational singularities. Then $\mathrm{t}_X \colon \mathbb{Q}_{\ell}(2)[4] \isoTo \mathscr{K}_X$ is an isomorphism in $\conDb{\etSite{X}, \mathbb{Q}_{\ell}}$. In particular, $\mathbb{Q}_{\ell}$ is a dualizing complex on $\etSite{X}$ in coefficient $\mathbb{Q}_{\ell}$.
\end{Cor}

\begin{Cor}[Poincar\'e Duality]
	Let $k$ be a separably closed field on which $\ell$ is invertible. And let $X$ be a surface over $k$. Assume that $X$ has only finitely many singular points, all of which are rational singularities. Then
	\begin{enumerate}
		\item The trace map $\Tr_X \colon \upHc{4}\bigl(X, \mathbb{Q}_{\ell}(2)\bigr) \isoTo \mathbb{Q}_{\ell}$ is an isomorphism.
		\item For every $0 \leqslant r \leqslant 4$ and every object $\mathscr{F}$ in $\conDb{\etSite{X}, \mathbb{Q}_{\ell}}$, the pairing
			\[ \upHc{r}(X, \mathscr{F}) \times \Ext^{4-r}_X\bigl(\mathscr{F}, \mathbb{Q}_{\ell}(2)\bigr) \to \upHc{4}\bigl(X, \mathbb{Q}_{\ell}(2)\bigr) \isoTo[\Tr_X] \mathbb{Q}_{\ell} \]
			is nondegenerate.
		\item For every $0 \leqslant r \leqslant 4$, the pairing
			\[ \upHc{r}(X, \mathbb{Q}_{\ell}) \times \upH{4-r}\bigl(X, \mathbb{Q}_{\ell}(2)\bigr) \to \upHc{4}\bigl(X, \mathbb{Q}_{\ell}(2)\bigr) \isoTo[\Tr_X] \mathbb{Q}_{\ell} \]
			is nondegenerate.
	\end{enumerate}
\end{Cor}

\begin{proof}
	(3) is a special case of (2). Let $p \colon X \to \Spec k$ be the structural morphism. By \autoref{M-A21} and \autoref{P-78E}, we have isomorphisms
	\[ \mathbb{Q}_{\ell}(2)[4] \isoTo[\mathrm{t}_X] \mathscr{K}_X \isoTo \rdF p^! \mathbb{Q}_{\ell} \,. \]
	Hence we have
	\begin{align*}
		\Ext_X^{4-r}\bigl(\mathscr{F}, \mathbb{Q}_{\ell}(2)\bigr) & \isoTo \Hom_{\conDb{\etSite{X}, \mathbb{Q}_{\ell}}}\bigl(\mathscr{F}[r], \mathbb{Q}_{\ell}(2)[4]\bigr) \\
		& \isoTo \Hom_{\conDb{\etSite{X}, \mathbb{Q}_{\ell}}}\bigl(\mathscr{F}[r], \rdF p^! \mathbb{Q}_{\ell}\bigr) \\
		& \isoTo \Hom_{\conDb{\mathbb{Q}_{\ell}}}\bigl(\rdF p_! \mathscr{F}[r], \mathbb{Q}_{\ell}\bigr) \\
		& \isoTo \upHc{r}(X, \mathscr{F})\spcheck \,.
	\end{align*}
	So (2) is proved. In above isomorphism, letting $r = 4$ and $\mathscr{F} = \mathbb{Q}_{\ell}$, we obtain
	\[ \upHc{4}(X, \mathbb{Q}_{\ell})\spcheck \isoTo \upH{0}\bigl(X, \mathbb{Q}_{\ell}(2)\bigr) = \mathbb{Q}_{\ell}(2) \,. \]
	So (1) is proved.
\end{proof}

\begin{Cor}
	Let $X$ be a surface having only finitely many singular points, all of which are rational singularities. Assume further that for any singular point $P$, $\mathcal{O}_{X, P}^{\mathrm{sh}}$ is factorial. Then $\mathrm{t}_X \colon \mathbb{Z}_{\ell}(2)[4] \isoTo \mathscr{K}_X$ is an isomorphism in $\conDb{\etSite{X}, \mathbb{Z}_{\ell}}$. In particular, $\mathbb{Z}_{\ell}$ is a dualizing complex on $\etSite{X}$ in coefficient $\mathbb{Z}_{\ell}$.
\end{Cor}

\begin{Rem}
	By \cite[(17.2)]{MR0276239}, the ring $\mathcal{O}_{X, P}^{\mathrm{sh}}$ is factorial if and only if its completion $\widehat{\mathcal{O}}_{X, P}^{\mathrm{sh}}$ is factorial. See \cite[\S25]{MR0276239} for the complete list of two-dimensional factorial complete local rings with rational singularities.
\end{Rem}

\section{Applications}

In this section, we study perverse sheaves and the weights of smooth sheaves on surfaces with rational singularities. First we briefly introduce the t-structures on arbitrary schemes.

Let $X$ be a scheme satisfying (\dag) equipped with a dimension function $\delta$. For each point $x$ on $X$, the canonical morphism $i_x \colon \Spec k(x) \to X$ factor as
\[ \Spec k(x) \xrightarrow{j} \overline{\{x\}} \xrightarrow{i} X \,. \]
Then we define two functors:
\begin{align*}
	i_x^{\ast} & \coloneqq j^{\ast} \circ i^{\ast} \colon \conDb{\etSite{X}, \overline{\mathbb{Q}}_{\ell}} \to \conDb(1){\etSite{x}, \overline{\mathbb{Q}}_{\ell}} \,, \\
	i_x^! & \coloneqq j^{\ast} \circ \rdF i^! \colon \conDb{\etSite{X}, \overline{\mathbb{Q}}_{\ell}} \to \conDb(1){\etSite{x}, \overline{\mathbb{Q}}_{\ell}} \,.
\end{align*}
Now we define a pair of full subcategories $\bigl(\lpD{\etSite{X}, \overline{\mathbb{Q}}_{\ell}}, \gpD{\etSite{X}, \overline{\mathbb{Q}}_{\ell}}\bigr)$ of $\conDb{\etSite{X}, \overline{\mathbb{Q}}_{\ell}}$ as following: for every object $\mathscr{F}$ in $\conDb{\etSite{X}, \overline{\mathbb{Q}}_{\ell}}$,
\begin{enumerate}
	\item $\mathscr{F} \in \lpD{\etSite{X}, \overline{\mathbb{Q}}_{\ell}}$ if and only if for any point $x$ on $X$, $i_x^{\ast}\mathscr{F} \in \conD[\leqslant -\delta(x)]{\etSite{x}, \overline{\mathbb{Q}}_{\ell}}$;
	\item $\mathscr{F} \in \gpD{\etSite{X}, \overline{\mathbb{Q}}_{\ell}}$ if and only if for any point $x$ on $X$, $i_x^!\mathscr{F} \in \conD[\geqslant -\delta(x)]{\etSite{x}, \overline{\mathbb{Q}}_{\ell}}$.
\end{enumerate}
By \cite{MR1971516} or \cite{MR2099084}, this pair defines a $t$-structure on $\conDb{\etSite{X}, \overline{\mathbb{Q}}_{\ell}}$. If $f \colon X \to Y$ is a quasi-finite morphism, then the functors $\rdF f_!$ and $f^{\ast}$ are right t-exact, and the functors $\rdF f^!$ and $\rdF f_{\ast}$ are left t-exact. Thus we may define the functor $f_{! \ast}$.

\begin{Prop}\label{D-P1}
	Let $X$ be a surface having only finitely many singular points, all of which are rational singularities; $U$ an open dense subset contained in $X_{\mathrm{reg}}$ and $j \colon U \hookrightarrow X$ the inclusion. Let $\mathscr{F}$ be a smooth $\overline{\mathbb{Q}}_{\ell}$-sheaf on $\etSite{X}$. Then $\mathscr{F}[2]$ is an object in $\mathrm{Perv}(X, \overline{\mathbb{Q}}_{\ell})$ and there is a canonical isomorphism $j_{! \ast} \circ j^{\ast}\mathscr{F}[2] \isoTo \mathscr{F}[2]$.
\end{Prop}

\begin{proof}
	Obviously we may assume that $U = X_{\mathrm{reg}}$. Then $Y \coloneqq X \setminus U$ contains only finitely many singular points, all of which are rational singularities of $X$. Let $P$ be a point in $X \setminus U$. Then
	\[ i_P^{\ast}\mathscr{F}[2] = \mathscr{F}_P[2] \in \lpD[-2]{\etSite{X}, \overline{\mathbb{Q}}_{\ell}} \,. \]
	Next we calculate $i_P^! \mathscr{F}[2]$. Put $\overline{X} \coloneqq \Spec \mathcal{O}_{X, P}^{\mathrm{sh}}$, and let $\rho \colon \overline{X} \to X$ and $\iota \colon P \hookrightarrow \overline{X}$ be the canonical morphisms. Since $\overline{X}$ is strictly local and $\mathscr{F}$ is a smooth sheaf, $\rho^{\ast}\mathscr{F} \isoTo \overline{\mathbb{Q}}_{\ell}^{\oplus r}$ for some $r \in \mathbb{N}$. So we have
	\begin{alignat*}{2}
		i_P^! \mathscr{F}[2] &= \rdF \iota^! \circ \rho^{\ast} \mathscr{F}[2] \isoTo \rdF \iota^! \overline{\mathbb{Q}}_{\ell}^{\oplus r}[2] \\
		& \isoTo \bigl(\rdF \iota^! \mathscr{K}_{\overline{X}}(-2)[-4]\bigr)^{\oplus r}[2] & & \qquad (\text{\autoref{M-A21}}) \\
		&= \rdF \iota^! \mathscr{K}_{\overline{X}}^{\oplus r}(-2)[-2] \\
		& \isoTo \mathscr{K}_P^{\oplus r}(-2)[-2] & & \qquad (\text{\autoref{P-78E}}) \\
		& \isoTo \overline{\mathbb{Q}}_{\ell}^{\oplus r}(-2)[-2] \,.
	\end{alignat*}
	Hence $i_P^! \mathscr{F}[2] \in \gpD[2]{\etSite{X}, \overline{\mathbb{Q}}_{\ell}}$. This show that $\mathscr{F}[2] \in \mathrm{Perv}(X, \overline{\mathbb{Q}}_{\ell})$.

	Let $i \colon Y \hookrightarrow X$ and $j \colon U \hookrightarrow X$ be the inclusion. Then $i^{\ast}\mathscr{F}[2] \in \lpD[-2]{\etSite{Y}, \overline{\mathbb{Q}}_{\ell}}$ and $\rdF i^! \mathscr{F}[2] \in \gpD[2]{\etSite{Y}, \overline{\mathbb{Q}}_{\ell}}$. So both ${}^{\mathrm{P}}i^{\ast}\mathscr{F}[2]$ and ${}^{\mathrm{P}}i^!\mathscr{F}[2]$ are zero. By \autoref{194L0}, there are two long exact sequences in $\mathrm{Perv}(X, \overline{\mathbb{Q}}_{\ell})$:
	\begin{gather*}
		\cdots \to {}^{\mathrm{P}}j_! \circ j^{\ast}\mathscr{F}[2] \to \mathscr{F}[2] \to i_{\ast} \circ {}^{\mathrm{P}}i^{\ast}\mathscr{F}[2] \to \cdots \,, \\
		\cdots \to i_{\ast} \circ {}^{\mathrm{P}}i^!\mathscr{F}[2] \to \mathscr{F}[2] \to {}^{\mathrm{P}}j_{\ast} \circ j^{\ast}\mathscr{F}[2] \to \cdots \,.
	\end{gather*}
	Thus ${}^{\mathrm{P}}j_! \circ j^{\ast}\mathscr{F}[2] \to \mathscr{F}[2]$ is an epimorphism and $\mathscr{F}[2] \to {}^{\mathrm{P}}j_{\ast} \circ j^{\ast}\mathscr{F}[2]$ is a monomorphism. Therefore $j_{! \ast} \circ j^{\ast}\mathscr{F}[2] = \mathscr{F}[2]$.
\end{proof}

Next we study the $\ell$-adic cohomology of surfaces over finite fields. Let $p \neq \ell$ be a prime number, $q \coloneqq p^e$ and $k \coloneqq \mathbb{F}_q$. If $X$ is an algebraic scheme over $k$, we write $\overline{X} \coloneqq X \otimes_k \bar{k}$; and for each sheaf $\mathscr{F}$ (or complex of sheaves) on $\etSite{X}$, we use $\overline{\mathscr{F}}$ to denote the inverse image of $\mathscr{F}$ on $\etSite{\overline{X}}$.

\begin{Defi}
	Let $w \in \mathbb{Z}$. An element $\alpha \in \overline{\mathbb{Q}}_{\ell}^{\ast}$ is said to be \emph{pure of weight $w$} if for any $\mathbb{Q}$-embedding $\iota \colon \overline{\mathbb{Q}}_{\ell} \hookrightarrow \mathbb{C}$, we have $|\iota(\alpha)| = q^{r/2}$.
\end{Defi}

\begin{Defi}
	Let $X$ be an algebraic scheme over $k$.
	\begin{enumerate}
		\item A constructible $\overline{\mathbb{Q}}_{\ell}$-sheaf $\mathscr{F}$ on $X$ is said to be
			\begin{itemize}
				\item \emph{punctually pure of weight $n$} if for every closed point $x$ on $X$, all eigenvalues of the geometric Frobenius $F_x$ acting on $\mathscr{F}_{\bar{x}}$ are pure of weight $w$;
				\item \emph{mixed} (\emph{of weights $\leqslant w$}), if there exists a finite filtration
					\[ 0 = \mathscr{F}_r \subset \mathscr{F}_{r-1} \subset \cdots \subset \mathscr{F}_1 \subset \mathscr{F}_0 = \mathscr{F} \]
					consisting of constructible $\overline{\mathbb{Q}}_{\ell}$-sheaves, such that each $\mathscr{F}_i/\mathscr{F}_{i+1}$ is punctually pure (of weights $\leqslant n$).
			\end{itemize}
		\item An object $\mathscr{F}$ in $\conDb{\etSite{X}, \widetilde{\mathbb{Q}}_{\ell}}$ is said to be
			\begin{itemize}
				\item \emph{mixed of weights $\leqslant w$} if for every $n \in \mathbb{Z}$, $\upH{n}(\mathscr{F})$ is mixed of weights $\leqslant n+w$;
				\item \emph{mixed of weights $\geqslant w$} if $\mathrm{D}(\mathscr{F}) = \mathbf{R}\mathscr{H}\!om(\mathscr{F}, \mathscr{K}_X)$ is mixed of weights $\leqslant -w$;
				\item \emph{pure of weight $w$} if it is mixed of weights $\leqslant w$ and $\geqslant w$.
			\end{itemize}
	\end{enumerate}
\end{Defi}

\begin{Prop}
	Let $X$ be a complete surface over $k$, $\mathscr{F}$ a smooth $\overline{\mathbb{Q}}_{\ell}$-sheaf on $\etSite{X}$. Assume that $X$ has only finitely many singular points, all of which are rational singularities; and $\mathscr{F}$ is punctually pure of weight $w$. Then for every $n \in \mathbb{N}$, the $\overline{\mathbb{Q}}_{\ell}$-module $\upH{n}(\etSite{\overline{X}}, \overline{\mathscr{F}})$ is punctually pure of weight $n+w$.
\end{Prop}

\begin{proof}
	Put $U \coloneqq X_{\mathrm{reg}}$ and let $j \colon U \hookrightarrow X$ be the inclusion. Since $U$ is smooth and $j^{\ast}\mathscr{F}$ is a smooth $\overline{\mathbb{Q}}_{\ell}$-sheaf, $j^{\ast}\mathscr{F}[2]$ is pure of weight $w+2$. By \autoref{D-P1}, $\mathscr{F}[2] \in \mathrm{Perv}(X, \overline{\mathbb{Q}}_{\ell})$ and $j_{! \ast} j^{\ast}\mathscr{F}[2] \cong \mathscr{F}[2]$. So after applying \cite[Corollaire 5.4.3]{MR751966}, we obtain that $j_{! \ast} j^{\ast}\mathscr{F}[2] \cong \mathscr{F}[2]$ is pure of weight $w+2$. Let $\pi \colon X \to \Spec k$ be the structural morphism. By \cite[(3.3.1) \& (6.2.3)]{MR601520}, the functor $\rdF \pi_{\ast} = \rdF \pi_!$ sends perverse sheaves pure of weight $w+2$ on $X$ to perverse sheaves pure of weight $w+2$ on $\Spec k$. Hence we prove the proposition.
\end{proof}

\bibliographystyle{plain}
\bibliography{dualizing}

\end{document}